\documentclass[12pt]{article}
\usepackage{amsmath,amssymb,latexsym}
\usepackage{fullpage}
\usepackage{theorem}

\begin{document}
\date{}
\title{\bf A  New Tower of Rankin-Selberg Integrals}
\author{ David Ginzburg\\ School of Mathematical Sciences\\
Sackler Faculty of Exact Sciences\\ Tel-Aviv University, Israel 69978\\
and\\ Joseph  Hundley \\ Penn State University Mathematics Department\\
University Park, State College PA, 16802, USA }
\maketitle
\baselineskip=18pt

\section{ Introduction}

The notion of a tower of Rankin-Selberg integrals was introduced
in \cite{G-R}. To recall this notion, let $G$ be a reductive group
defined over a global field $F$. Let $^LG$ denote the $L$ group of
$G$. Let $\rho$ denote a finite dimensional irreducible
representation of $^LG$. Given an irreducible generic cuspidal
representation of $G({\bf A})$, we let $L^S(\pi,\rho,s)$ denote
the partial $L$ function associated with $\pi$ and $\rho$. Here
$s$ is a complex variable and ${\bf A}$ denotes the adele ring
associated with $F$. If $\rho$ acts on the vector space $V$, we
denote by ${\bf C}[V]$ the symmetric algebra attached to the
vector space $V$. Let ${\bf C}[V]^{^LG}$ denote the $^LG$
invariant polynomials inside the symmetric algebra. As far as we
know all examples of $L$ functions represented by a Rankin-Selberg
integral are associated with representations $\rho$ such that
${\bf C}[V]^{^LG}$ is a free algebra. A list of all such groups,
representations and the degrees of the generators of the invariant
polynomials are given in \cite{K}.

The basic observation in \cite{G-R} is that there is some relation
between the Eisenstein series one uses to construct the
Rankin-Selberg integral and the number of generators of the
invariant polynomials and their degrees. This relation is far from
being clear and it is mainly based on observation of all known
constructions of such integrals. To summarize in an unprecise
manner, the relations
are:\\
{\bf 1)}  If $\rho_1$ and $\rho_2$ have the same number of
generators with the same degrees, then in some cases the
Rankin-Selberg integrals which represent the corresponding two
$L$ functions, use the same Eisenstein series.\\
{\bf 2)} Suppose that the Eisenstein series one uses for a certain
construction is defined over  $H({\bf A})$, where $H$ is a
reductive group. Suppose that this Eisenstein series corresponds
to an induced  representation induced from a parabolic subgroup
$P=MU$ of $H$. Here $M$ is the Levi part of $P$ and $U$ its
unipotent radical. The group $M({\bf C})$ acts on $U({\bf C})$ by
conjugation and one obtains this way $r$ irreducible finite
dimensional representations of $M({\bf C})$. Suppose that the
corresponding Rankin-Selberg integral represents the $L$ function
$L^S(\pi,\rho,s)$. Let $k$ denote the number of generators of
${\bf C}[V]^{^LG}$, where we recall that we assume that ${\bf
C}[V]^{^LG}$ is free. Then the second observation that was pointed
out in \cite{G-R} is that $r\ge k$.

It should be clear that these two observations are based mainly on
experience and we are not aware of precise theoretical reasons. We
also want to mention that information on $L$ functions
$L^S(\pi,\rho,s)$ where $\rho$ does not satisfy the above
properties, can be obtained using other methods such as lifting
theory.

In this paper we wish to point out two more observations that may
shed some more light on the above relations. It will be convenient
to first illustrate these observations by considering two
examples.

Consider the following example of a tower given in \cite{G-R}:
$$\begin{matrix} ({\bf a1})&&&G=GL_{n}&&& ^LG=GL_{n}({\bf C})&&&
\rho=2\varpi_1\\
({\bf a2})&&&G=GL_{n}\times GL_n&&& ^LG=GL_{n}({\bf C}) \times
GL_{n}({\bf C})&&& \rho=\varpi_1\times \varpi_1\\
({\bf a3})&&&G=GL_{2n}&&& ^LG=GL_{2n}({\bf C})&&& \rho=\varpi_2
\end{matrix}$$

We recall the construction of the Rankin-Selberg integral which
represents the $L$ function in case $({\bf a3})$. This integral
was introduced in \cite{J-S} and is given by
\begin{equation}\label{int1}
\int\limits_{Z({\bf A})GL_n(F)\backslash GL_n({\bf A})}
\int\limits_{Mat_{n\times n}(F)\backslash Mat_{n\times n}({\bf
A})}\varphi_\pi\left (
\begin{pmatrix} I&X\\&I \end{pmatrix} \begin{pmatrix} g&\\&g
\end{pmatrix}\right ) E(g,s)\psi(trX)dXdg
\end{equation}
Here $\varphi_\pi$ is a vector in the space of $\pi$ which is an
irreducible cuspidal representation defined on $GL_{2n}({\bf A})$,
and $E(g,s)$ is an Eisenstein series defined on the group
$GL_n({\bf A})$. For more details see \cite{J-S}. Let us show how
the integral which represents the $L$ function given in $({\bf
a2})$ can be derived from integral \eqref{int1}. First notice that
$GL_n\times GL_n$ is a Levi part of a maximal parabolic subgroup
$P$ of $GL_{2n}$. Now suppose we formally replace in \eqref{int1}
the cuspidal representation $\pi$ by the Eisenstein series
$E_{\tau,\sigma}(g,\nu)$ associated with the induced
representation $Ind_{P({\bf A})}^{GL_{2n}({\bf A})}(\tau\otimes
\sigma)\delta_P^\nu.$ Here $\tau$ and $\sigma$ are cuspidal
representations defined on $GL_n({\bf A})$ and $\nu$ is a complex
variable. Of course the integral will not converge. However, if we
ignore this issue, and formally unfold the Eisenstein series
$E_{\tau,\sigma}(g,\nu)$, we are led to consider the space of
double cosets $P\backslash GL_{2n}/GL_nX$. If we consider the open
orbit contribution to the integral, it is not hard to check that
we obtain the integral
\begin{equation}\label{int2}
\int\limits_{Z({\bf A})GL_n(F)\backslash GL_n({\bf A})}
\varphi_\tau(g)\varphi_\sigma(g)E(g,s)dg
\end{equation}
as inner integration. As is well known, integral \eqref{int2}
represents the tensor product $L$ function of $\tau\times\sigma$.
In other words, this integral is the one which represents the $L$
function described in case $({\bf a2})$. Furthermore, if one
restricts the exterior square representation $\varpi_2'$ of
$GL_{2n}({\bf C})$ to $GL_n({\bf C})\times GL_n({\bf C})$, then
one obtains $\varpi_2'|_{GL_n\times
GL_n}=(\varpi_1\times\varpi_1)\oplus(\varpi_2\times 1)\oplus
(1\times\varpi_2)$. From this we deduce the following. If we start
with the representation $\rho$ as defined in case $({\bf a3})$ and
restrict it to the $L$ group of the Levi part then the
representation $\rho$ corresponding to case $({\bf a2})$ occurs in
the restriction. Moreover its the representation with the largest
dimension which occurs in the restriction.

The formal replacement of a cuspidal representation by an
Eisenstein series and then analyzing the contribution from the
open orbit is one of the observations we wish to make. It should be
mentioned that this observation does not explain how to derive a
global construction that will represent the $L$ function described
in case $({\bf a1)}$. We now consider the second example of a
tower as described in \cite{G-R}. This tower consists of four
members as follows:
$$\begin{matrix} ({\bf b1})&&& G=GL_2&&& ^LG=GL_2({\bf C})&&&
\rho=4\varpi_1\\
({\bf b2})&&& G=GL_3&&& ^LG=GL_3({\bf C})&&&
\rho=\varpi_1+\varpi_2\\
({\bf b3})&&& G=GSpin_7&&& ^LG=GSp_6({\bf C})&&& \rho=\varpi_2\\
({\bf b4})&&& G=F_4&&& ^LG=F_4({\bf C})&&& \rho=\varpi_4
\end{matrix}$$

The construction of Rankin-Selberg integrals for cases $({\bf
b1}),({\bf b3})$ and $({\bf b4})$ was given in \cite{G-R}. The
case $({\bf b2})$ was studied in \cite{G1}. The integral which
represents the $L$ function given in $ ({\bf b4})$ can be
described as follows. Let $\pi$ denote a generic cuspidal
representation defined on the group $F_4({\bf A})$. Let $E(g,s)$
denote the degenerate Eisenstein series defined on the exceptional
group $G_2({\bf A})$ as described in section 1 in \cite{G-R}. The
global integral is
\begin{equation}\label{int3}
\int\limits_{G_2(F)\backslash G_2({\bf
A})}\int\limits_{U(F)\backslash U({\bf
A})}\varphi_\pi(ug)E(g,s)\psi_U(u)dudg
\end{equation}
Here $U$ is a certain unipotent subgroup of $F_4$ and $\psi_U$ is
an additive character defined on the group $U$. Observe that
$GSpin_7$ is a Levi part of a maximal parabolic subgroup $P$ of
$F_4$. Let $\tau$ denote a generic cuspidal representation defined
on the group $GSpin_7$. Let $E_\tau(g,\nu)$ denote the Eisenstein
series defined on the group $F_4({\bf A})$ which is associated to
the induced representation $Ind_{P({\bf A})}^{F_{4}({\bf
A})}\tau\delta_P^\nu.$ If we formally replace in \eqref{int3} the
cuspidal representation $\pi$ by $E_\tau(g,\nu)$, and then unfold
this Eisenstein series, then we obtain from the open orbit
\begin{equation}\label{int4}
\int\limits_{G_2(F)\backslash G_2({\bf A})}
\varphi_\tau(g)E(g,s)dg
\end{equation}
as inner integration. As described in \cite{G-R} section 4 this is
precisely the global integral which represents the $L$ function
which is described in $({\bf b3)}$. Further more, let $Q$ denote
the maximal parabolic subgroup of $Spin_7$ whose Levi part is
$GL_3$. Let $\sigma$ denote a cuspidal representation defined on
the group $GL_3({\bf A})$. Replace in \eqref{int4} the cuspidal
representation $\tau$ by the Eisenstein series $E_\sigma(g,\nu)$
which is associated with the induced representation $Ind_{Q({\bf
A})}^{Spin_7({\bf A})}\sigma\delta_Q^\nu.$ Unfolding the integral
we obtain from the open orbit
\begin{equation}\label{int5}
\int\limits_{SL_3(F)\backslash SL_3({\bf A})}
\varphi_\sigma(g)E(g,s)dg
\end{equation}
as inner integration. As described in \cite{G1} this is precisely
the global integral which represents the $L$ function described in
$({\bf b2})$.

As in the previous case we can restrict in each case the
representations $\rho$ to the $L$ group of the Levi part. Suppose
that $(0,0,0,1)$ is the representation of $F_4({\bf C})$ of
dimension 26. This is the representation $\rho$ obtained in case
$({\bf b4})$. Restrict it to $GSp_6({\bf C})$ which is the $L$
group of $GSpin_7$. We obtain $(0,0,0,1)|_{C_3}=(0,1,0)+2(1,0,0)$.
Here $(0,1,0)$ is the second fundamental representation of
$GSp_6({\bf C})$ which has degree 14, and $(1,0,0)$ is the six
dimensional standard representation. If we further restrict
$GSp_6({\bf C})$ to $GL_3({\bf C})$ we obtain
$(0,1,0)|_{A_2}=(1,1)+(1,0)+(0,1)$. Again, as in the first tower
we can see that if we restrict $\rho$ as defined in case $({\bf
b4})$ we obtain the representation $\rho$ as defined in case
$({\bf b3})$ as the largest piece in the restriction. Similarly,
if we restrict from case $({\bf b3})$ to $({\bf b2})$.

We mention again that this observation does not allow one to
obtain the integrals for cases $({\bf a1})$ and 
$({\bf b1})$.  The construction in these cases is more 
complicated and involves covering groups.  

To summarize, the above examples suggests the following two
points:\\
{\bf 1)} Suppose that we are given a Rankin-Selberg integral which
we know how to unfold to an Eulerian integral with the Whittaker
function defined on the cuspidal representations. Then replacing a
cuspidal representation by an Eisenstein series and considering
the contribution from the open orbit, sometimes yields a new
Eulerian Rankin-Selberg integral. In fact, one can replace the
cuspidal representation by various Eisenstein series. Experience
indicates that most of the time one gets either zero, or an
integral which does not unfold to a Whittaker integral. The second
point is\\
{\bf 2)} Suppose that the Eulerian integral we start with
represents an $L$ function which is associated to the finite
dimensional irreducible representation $\rho$ of the complex group
$^LG$. Suppose that we replace a cuspidal representation, defined
over the group $G({\bf A})$, by an Eisenstein series induced from
a cuspidal representation defined on the Levi part $M({\bf A})$.
Suppose that when we formally unfold the new integral, the
contribution from the open orbit produces a new integral which is
Eulerian with Whittaker functions. Then the new integral will
represent the $L$ function associated with the largest irreducible
representation which occurs in the restriction $\rho|_{^LM}$.

In these notes we announce a construction of a new tower of
Rankin-Selberg integrals. The tower we consider is the following

$$\begin{matrix} ({\bf c1})& G=GL_3\times GL_2& ^LG=
GL_3({\bf C})\times GL_2({\bf C})& \rho=2\varpi_1\times
\varpi_1\\
({\bf c2})& G=GL_3\times GL_3\times GL_2& ^LG=GL_3({\bf C})\times
GL_3({\bf C})\times GL_2({\bf C})&
\rho=\varpi_1\times\varpi_1\times\varpi_1\\
({\bf c3})& G=GL_6\times GL_2& ^LG= GL_6({\bf C})\times
GL_2({\bf C})& \rho=\varpi_2\times\varpi_1\\
({\bf c4})& G=E_6\times GL_2& ^LG= E_6({\bf C})\times GL_2({\bf
C})& \rho=\varpi_1\times\varpi_1 \end{matrix}$$

It follows from \cite{K} that in all these representations the
$^LG$ invariant algebra has one generator of degree 12. At this
point we know a Rankin-Selberg construction for all three cases
$({\bf c2})-({\bf c4})$. In the next section we shall explain
these constructions and show in an example  how to derive one
integral from the other. One can also check that restricting from
one case to the other does indeed produce the right representation
$\rho$ in each case.

It should be mentioned that all of these $L$-functions can be
studied using the Langlands-Shahidi method as explained in
\cite{S}.

\section{ The Global Integrals}

We start with the global construction which will correspond to
case $({\bf c4})$, as explained in the introduction. Let $G$
denote the similitude exceptional group of type $E_6$,
constructed exactly as in \cite{G2}. To
introduce the global integral we shall need to consider two small
representations which we shall now define. First, let $\theta$
denote the minimal representation defined on $G({\bf A})$. This
representation was constructed and studied in \cite{G-R-S}. The
construction there is defined on the group $E_6$, however there
are no problems to extend this definition to similitude groups.
See \cite{G-J} for a similar definition for the similitude
exceptional group $GE_7$. In this paper we shall denote a function
in the space of this representation by $\theta(g)$. Another
representation we will need for our construction was defined and
studied in \cite{G-H} section 3. The representation constructed
there was defined on the group $GSO_{10}({\bf A}).$ A similar
definition holds for the group $GSpin_{10}({\bf A}).$ This
representation depends on a cuspidal representation $\tau$ defined
on $GL_2({\bf A}).$ We shall denote a vector in this space by
$\theta_\tau(h)$ where $h\in GSpin_{10}({\bf A}).$ We briefly
recall the definition. Let $R$ denote the parabolic subgroup of
$GSpin_{10}$ whose Levi part is $GL_3\times GSpin_4$. Let
$\epsilon(\tau)=\tau\otimes\tau$ and let $\mu(\tau)$ denote the
symmetric square lift of $\tau$ to $GL_3$ as constructed in
\cite{Ge-J}. Let $E(\tau,h,s)$ denote the Eisenstein series
defined on $GSpin_{10}({\bf A})$ associated with the induced
representation $Ind_{R({\bf A})}^{GSpin_{10}({\bf
A})}(\mu(\tau)\otimes\epsilon(\tau))\delta_R^s$. It is not hard to
check that this Eisenstein series has a unique simple pole, and we
denote the residue representation by $\theta_\tau$.

Using this last representation, we shall now construct the
Eisenstein series we use in our global construction. Let $P$
denote the maximal standard parabolic subgroup of $G$ whose 
Levi part contains all the simple roots except $\alpha_1$.  
This Levi part is essentially $GSpin_{10}$.
Let $E_\tau(g,s)$ denote the Eisenstein series
defined on $G({\bf A})$ which is associated to the induced
representation $Ind_{P({\bf A})}^{G({\bf
A})}\theta_\tau\delta_P^s$.

Let $\pi$ denote a generic cuspidal representation defined on
$G({\bf A})$. We shall assume that $\pi$ has a trivial central
character. Consider the global integral
\begin{equation}\label{glob1}
\int\limits_{Z({\bf A})G(F)\backslash G({\bf A})}\varphi_\pi(g)
\theta(g)E_\tau(g,s)dg
\end{equation}
Here $Z$ denotes the center of $G$ and $\varphi_\pi$ is a vector
in the space of $\pi$. This integral represents the $L$ function
corresponding to case $({\bf c4})$.

Let us show how to obtain the Rankin-Selberg integral which will
represent case $({\bf c3})$ as denoted in the introduction. Let
$Q$ denote the maximal parabolic subgroup of $G$ whose Levi part
is $M=GL_1\times GL_6$. Let $\sigma$ denote a cuspidal
representation of $GL_6({\bf A})$ with trivial central character.
Let $E_\sigma(g,\nu)$ denote the Eisenstein series defined on
$G({\bf A})$ associated with the induced representation
$Ind_{Q({\bf A})}^{G({\bf A})}\sigma \delta_Q^\nu$. In
\eqref{glob1} we replace the function $\varphi_\pi(g)$ by
$E_\sigma(g,\nu)$. Even though the integral does not converge, we
formally unfold the Eisenstein series $E_\sigma(g,\nu)$ to obtain
\begin{equation}\label{glob2}
\int\limits_{Z({\bf A})M(F)U(F)\backslash G({\bf
A})}f_\sigma(g,\nu) \theta(g)E_\tau(g,s)dg
\end{equation}
Here $U$ is the unipotent radical of $Q$ and $f_\sigma(g,\nu)$
defines a section in the corresponding induced representation.
Recall that $U$ has a structure of a Heisenberg group with 21
variables. Let $x_{122321}(r)$ denote the one dimensional
unipotent subgroup $U$ which is the center of $U$. Here, and
henceforth we shall use the notations for various roots of the
group $G$ as defined in \cite{G2}. We expand $\theta(g)$ along the
center of $U$. That is, we expand it along the unipotent group
generated by $x_{122321}(r)$ with points in $F\backslash {\bf A}$.
The group $M(F)$ acts on this expansion with two orbits. Ignoring
the trivial orbit, we obtain the contribution
\begin{equation}\label{glob3}
\int\limits_{Z({\bf A})H(F)U(F)\backslash G({\bf
A})}\int\limits_{F\backslash {\bf A}}
\theta(x_{122321}(r_1)g)\psi(r_1)dr_1f_\sigma(g,\nu)E_\tau(g,s)dg
\end{equation}
Here $H$ is the stabilizer inside $M$ of the character $\psi$. One
can check that $H=\{g\in GL_6 : detg\ \ \text{is a square}\}$.
Factoring the integration over $H$ and over the center of $U$ we
obtain after a change of variables, the integral
\begin{equation}\label{glob4}
\int\limits_{Z({\bf A})H(F)\backslash H({\bf
A})}\int\limits_{U(F)\backslash U({\bf
A})}\int\limits_{(F\backslash {\bf A})^2}
\varphi_\sigma(h)\theta(ux_{122321}(r_1)h)
E_\tau(ux_{122321}(r_2)h,s)\psi(r_1-r_2)dr_1dr_2dh
\end{equation}
as inner integration. Here $\varphi_\sigma$ is a vector in the
space of the cuspidal representation $\sigma$. Notice that this
integral converges absolutely. This is our candidate for the
Rankin-Selberg integral which will represent case $({\bf c3})$.

We can further continue and replace $\sigma$ by an Eisenstein
series. Indeed let $\pi_1$ and $\pi_2$ denote two cuspidal
representations of $GL_3({\bf A})$. Let $L$ denote the parabolic
subgroup of $GL_6$ whose Levi part is $GL_3\times GL_3$. Let
$E_{\pi_1,\pi_2}(x,\nu)$ denote the Eisenstein series associated
with the induced representation $Ind_{L({\bf A})}^{GL_6({\bf
A})}(\pi_1\otimes\pi_2) \delta_L^\nu$. Replacing in \eqref{glob4}
the cuspidal representation $\sigma$ by this Eisenstein series (
again, this is a formal process, since the integral does not
converge) and performing certain Fourier expansions, one obtains
the integral
\begin{align}\label{glob5}
\int\limits_{Z({\bf A})H(F)\backslash H({\bf A})}\int\limits_{
V(F)\backslash V({\bf A})}\int\limits_{(F\backslash {\bf
A})^3}\varphi_{\pi_1,\pi_2}(h)\theta(x_{010000}(r_1)vx_{112321}(r_2)&x_{122321}
(r_3)h)\times \\
E_\tau(x_{010000}(r_1)vh,s)\psi(r_1+r_2)dr_idvdh\notag
\end{align}
as inner integration. Here $\varphi_{\pi_1,\pi_2}$ is a vector in
the space of $\pi_1\otimes\pi_2$. We also have $H=\{(g_1,g_2)\in
GL_3\times GL_3 : detg_1=detg_2\}$ and the group $V$ is the
standard  unipotent radical of the maximal parabolic subgroup of
$G$ whose Levi part is $GL_3\times GL_3\times GL_2$. This is the
integral which represents the case $({\bf c2})$.

At this point we unfolded all these three integrals and
established that they are indeed Eulerian. This we achieved by
obtaining the Whittaker function of each of the cuspidal
representations involved in the integral. As always, with these
type of integrals, the unfolding process is long and tedious but
quite straightforward. The next step is to compute the unramified
local integrals. So far we have performed some of the calculations
which indicate that our integrals do represent the $L$ functions
in question. It is not yet clear to us how complicated will be the
decomposition of the symmetric algebras in the various cases. It
will also be interesting to study the possible poles of these $L$
functions. This will be accomplished by understanding the poles of
the Eisenstein series we use in all these cases.

We are also interested in finding the Rankin-Selberg integral
which represents case $({\bf c1})$. Past experience indicates that
some of the representations involved should be defined on a
covering group. So far we don't know how to do it.

We summarize\\
{\bf Theorem:} {\sl Integrals \eqref{glob1}, \eqref{glob4} and
\eqref{glob5} are Eulerian. Each of these three integrals unfolds
to the Whittaker function defined on each cuspidal representation
which appears in the integral. Integral \eqref{glob1} represents
the partial $L$ function $L^S(\pi\times\tau,St\times St,s)$ where
$St\times St$ corresponds to the standard representation of
$^LG\times GL_2({\bf C})$. Integral \eqref{glob4} represents
$L^S(\sigma\times\tau,\wedge^2\times St,s)$ where $\wedge^2$ is
the exterior square representation of $GL_6({\bf C})$, and
integral \eqref{glob5} represents $L^S(\pi_1\times\pi_2
\times\tau,St\times St\times St,s)$.}

\end{document}